\documentclass[11pt]{amsart}

\usepackage{amssymb, amsfonts, amscd}
\usepackage{mathtools} % includes amsmath

\usepackage{graphicx, epsf, svg, standalone, tikz}
\usepackage{wrapfig}
\usepackage{adjustbox}
\usepackage[table]{xcolor}

\usepackage{parskip}
\usepackage{enumitem}
\usepackage{float}
\usepackage{caption}
\usepackage{comment}
\usepackage{threeparttable}
\usepackage{arcs}

\usepackage[most]{tcolorbox} % modern replacement for mdframed

\tcbset{
  exampledefault/.style={
    colback=white,
    colframe=black,
    boxrule=0.8pt,
    left=5mm,
    right=5mm,
    top=2mm,
    bottom=2mm,
    arc=0mm,
    outer arc=0mm,
    boxsep=3pt,
    breakable,
    enhanced,
  },
}

\usepackage{hyperref}
\hypersetup{
  colorlinks=true,
  linkcolor=blue,
  filecolor=magenta,
  urlcolor=red,
  citecolor=gray,
}

\usepackage[capitalise]{cleveref}

\definecolor{green}{RGB}{50,205,50}

\begin{document}

\newtheorem{defi}{Definition}
\newtheorem{thm}{Theorem}
\newtheorem*{conjecture}{Conjecture}
\newtheorem{prop}{Proposition}
\newtheorem{lemma}{Lemma}
\newtheorem{coro}{Corollary}
\newtheorem{quest}{Question}
\newtheorem{corollary}{Corollary}
\newtheorem{problem}{Problem}
\newtheorem{ex}{Example}
\newtheorem*{theorem*}{Theorem}

% Define a new theorem style for remarks
\newtheoremstyle{boldremark}
  {\topsep}   % Space above
  {\topsep}   % Space below
  {}          % Body font (default is upright)
  {}          % Indent amount
  {\bfseries} % Theorem head font (e.g., Remark in bold)
  {.}         % Punctuation after theorem head
  { }         % Space after theorem head
  {}          % Theorem head spec

% Apply the new theorem style to the remark environments
\theoremstyle{boldremark}
\newtheorem*{rem}{Remark}
\newtheorem{numrem}{Remark}

\newtheorem{claim}{Claim}
\newenvironment{proofc}{\begin{proof}[Proof of Claim]\renewcommand*{\qedsymbol}{\(\boxtimes\)}}{\end{proof}}

\newcommand{\claimqed}{\hfill $\blacksquare$}

\title{Shrinking the Jung radius: Maximizing partial coverage of finite point sets}
\author{Andr\'{a}s Bezdek}
\address{Andr\'{a}s Bezdek, Mathematics \& Statistics, Auburn Univ., Auburn, AL 36849-5310, USA}
\email{bezdean@auburn.edu}

\author{Owen Henderschedt}
\address{Owen Henderschedt, Mathematics \& Statistics, Auburn Univ., Auburn, AL 36849-5310, USA}
\email{olh0011@auburn.edu}

% \centerline{\today}


\begin{abstract}
Jung's theorem says that planar sets of diameter $1$ can be covered by a closed circular disk of radius $\frac 1{\sqrt3}$. In this paper we  consider a fractional Jung-type problem for finite planar point-sets. Let  $\mathcal{P}_n$ be the family of all finite sets of $n$ points in the plane, of diameter at most $1$. Let the function value $N_n(r)$ ($0 < r \leq 1$) be the largest integer $k$ so that for every point set $P \in  \mathcal{P}_n$ there is a closed circular disk of radius $r$ which covers at least $k$ points of $P$. We focus on the radii $r=\frac 12$ and $r=\frac 14$ and prove exact maximum values. Concerning the radius $r= \frac 12$, we prove $N_n(\frac{1}{2})=\lceil \frac{n}{3}\rceil+1$. Concerning the radius $r= \frac 14$, we prove that 
$N_{n}(\frac{1}{4}) = \lceil \frac{n}{7}\rceil$ if $n$ is not a multiple of 7, and  $N_{n}(\frac{1}{4})$ is $ \frac{n}{7}$ or $ \frac{n}{7}+1$ otherwise. We also initiate further study of the function $N_n(r)$ by giving  lower and upper bounds for $N_n(r)$  ($0 < r < \frac 1{\sqrt3}$). 
\end{abstract}

\maketitle

\section{Introduction}\label{sec:intro}

\noindent At the turn of the 20th century, Heinrich Jung proved the following celebrated theorem:

\begin{thm}\label{jung}(Jung's theorem in the plane \cite{HJ01, HJ10}) Given any set of points $P$ in the plane such that the distance between any two points is at most $1$, there exists a circular disk (in short a circle) with radius $\frac{1}{\sqrt{3}}$ that can cover every point in $P$. \end{thm}

Since $\frac{1}{\sqrt{3}}$ is the circumradius of an equilateral triangle with side length $1$, no circle of radius smaller than $\frac{1}{\sqrt{3}}$ can cover all three vertices of the triangle. This means that if the radius $\frac{1}{\sqrt{3}}$ in Jung's theorem shrinks, it is very natural to ask what fraction of the points can always be covered by a circle of radius $r$. More formally:

\begin{defi}
Let $n$ be a fixed positive integer. Let 
$\mathcal{P}_n$ be the family of all sets of $n$ points, such that in any set, the distance between any two points is at most $1$. We define the function value  $N_n(r)$, for $0 < r \leq 1$, to be the largest integer $k$ so that for every point set $P$  $\in \mathcal{P}_n$ there exists a circle of radius $r$ which  covers at least $k$ points in $P$. 
\end{defi}

Using this notation, Jung's Theorem says that  $N_n(r) = n $ for any $\frac 1{\sqrt{3}} \leq r\leq 1$. Our goal in this paper is to begin investigating the values of $N_n(r)$ when $r<\frac 1{\sqrt{3}}$. Determining $N_n(r)$ for all possible values of $r$ appears to be quite difficult. Thus, it is natural to ask: Which radii $r$ are good candidates for determining the exact function values $N_n(r)$? Both $r=\frac{1}{2}$ and $r=\frac{1}{4}$ seem to be natural choices, in view of some previously related results.

We start with $r = \frac{1}{2}$. Problem 1.20 (due to K.~Vasilev and K.~Garov) of the 2023 Bulgarian Math Olympiad \cite{BMO23} essentially asks students to prove that for sufficiently small $\varepsilon > 0$, $N_{3n}(r) = n$ whenever $r \in \left(\frac{1}{2} - \varepsilon, \frac{1}{2}\right)$. Interestingly, the case $r = \frac{1}{2}$ itself is not included in Problem~1.20, and for good reason. Not only is it the more difficult case, but it also yields a different outcome: as it turns out, $N_{3n}\!\left(\frac{1}{2}\right) = n + 1$. This follows as a corollary of one of the main results presented in this paper.

\begin{thm}\label{thm:half}
Given any set of $n$ points in the plane such that the distance between any two points is at most $1$, there exists a circle of radius $\frac{1}{2}$ that covers at least $\left\lceil \frac{n}{3} \right\rceil + 1$ of the points. Moreover, this bound is best possible. That is, $N_n\left(\frac{1}{2}\right) = \left\lceil \frac{n}{3} \right\rceil + 1$.
\end{thm}

Moreover, using an approach different from the posted solution in \cite{BMO23}, we are able to extend this result and show that $N_{3n}(r) = n$ holds for the larger interval $\left(\frac{\sqrt{3}}{4}, \frac{1}{2}\right)$ (see Theorem~\ref{thm: B}).

Next, we  describe a general idea for finding lower and upper bounds for $N_n(r)$ for a selected set of radii $r$. It will turn out that $r= \frac 14$ is among these radii, which explains why $r = \frac14$ became our second choice for finding exact values of $N_n(r)$.

\begin{figure}[h!]
\begin{center}
\includegraphics[scale=.75]{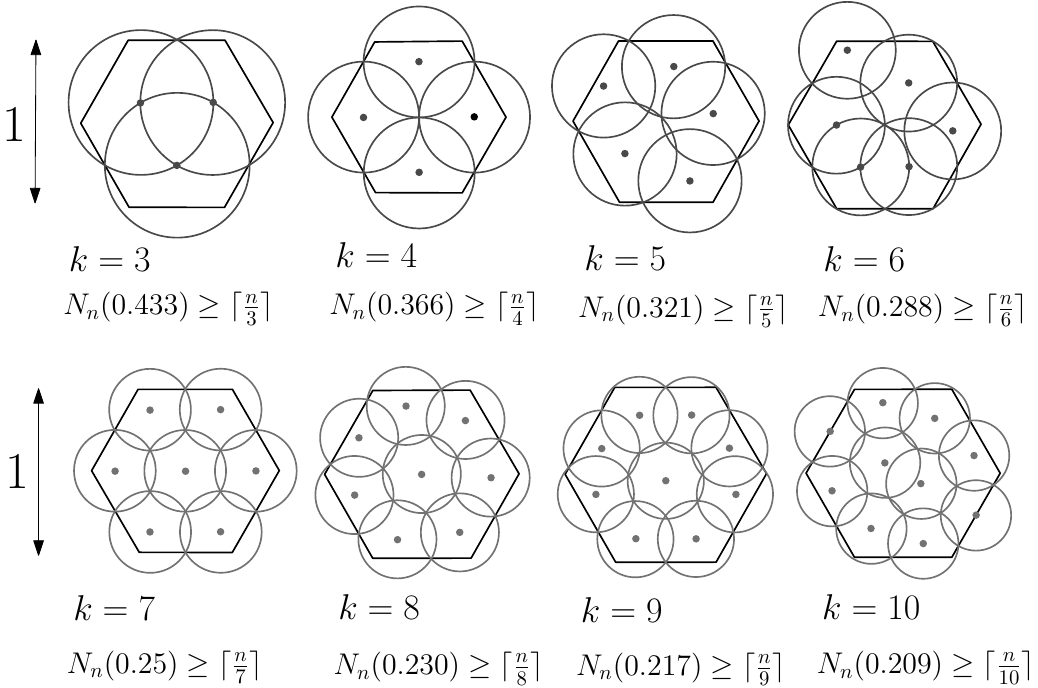}
\caption{Circle coverings of a regular hexagon found by Y. Liu in \cite{Y22}}
\label{Yanchao}
\end{center}
\end{figure}

\noindent {\bf How to establish lower bounds for $N_n(r)$?} We start with a convex disk, capable of covering every point set of diameter one, otherwise known as a {\it universal cover}. It is well known that the regular hexagon whose opposite sides are at unit distance is a universal cover. Similarly, the square of edge length 1 is also a  universal cover. If we can then cover a universal cover with $k$ congruent circles of radius $r$, then by applying the pigeonhole principle we obtain that one of these circles contains at least $\lceil \frac nk \rceil$ points. Thus, we readily have $N_n(r) \geq \lceil \frac nk \rceil$. Fortunately, numerous papers deal with such finite circle coverings. For example, the circle coverings in Figure \ref{Yanchao} come from a paper by Y. Liu \cite{Y22} whose goal was to find the smallest radius such that $k$ congruent circles can cover a specific regular hexagon. Liu approached this problem algorithmically and obtained improved bounds for $k$ between $3$ and $10$. By scaling the hexagons appropriately, we also include in Figure \ref{Yanchao} for each circle covering   the number of circles $k$ in the covering and  the approximate decimal value of the radii $r$ where pigeonhole principle gives  the lower bound of $N_n(r_k)\geq \lceil \frac nk \rceil$. Certainly, $\lceil \frac nk \rceil$ is a lower bound for $N_n(r)$ for all $r \geq r_k$.

Looking at the circle covering of seven circles (i.e. $k=7$ in Figure \ref{Yanchao}), we get the lower bound $\lceil \frac{n}{7}\rceil$ when $r = \frac14$. We will show that this bound is also an upper bound for point sets whose number of points is not a multiple of $7$. In fact, we show more.

\begin{thm}\label{thm:fourth}
$N_{n}(\frac{1}{4}) = \lceil \frac{n}{7}\rceil$ if $n$ is not a multiple of 7, and  $N_{n}(\frac{1}{4})$ is $ \frac{n}{7}$ or $ \frac{n}{7}+1$ otherwise. 
\end{thm}

 For completeness, we mention that 7 is the only multiple of seven where it is known which of $ \frac{n}{7}$ and $ \frac{n}{7}+1$ is the actual  maximum.  P. Bateman and P. Erd\H{o}s \cite{BE51} showed that $N_7(\frac{1}{4})=\frac 77 +1=2$ (see \cref{sec:fourth} for more on this).

\noindent {\bf How to establish upper bounds for $N_n(r)$?} To show $N_n(r)\leq k$, all we need is a concrete arrangement of $n$ points so that no circle of radius $r$ can cover more than $k$ points. The search for  special point sets is at the heart of our paper.  We will refer to, and label, these point sets as {\it examples}. To help the reader navigate these constructions, each example is placed in a text box for clarity and visual separation.

In several of our examples, it is common to construct an $n$-point set by taking an $s$-point set $S$, for some $s < n$, and then duplicating or stacking each point some number of times. The issue with this approach is that $n$ may not be an integer multiple of $s$. If this is the case, we duplicate the entire point set $\lfloor \frac{n}{s} \rfloor$ times, and then, for the remaining $n - \lfloor \frac{n}{s} \rfloor s$ points, we place each one on a unique existing stack. Hence, the number of points in any two stacks differs by at most $1$. We call this process constructing a \emph{balanced $S$}. Note that if $N_s(r) \leq c$, then this implies $N_n(r) \leq \lceil \frac{nc}{s} \rceil$. 

\noindent Due to the stacking nature of these constructions, for a fixed $n$, values of $N_n(r)$ are influenced by various divisibility properties of $n$ which explains the use of ceiling functions in our theorems. To complement this, we also consider the asymptotic behavior of $c(r) := \lim_{n\to \infty} \frac{N_n(r)}{n}$. Our final contribution to the study of this fractional version of Jung’s theorem is consolidated in the following result, which packages together all currently known bounds and exact values of the asymptotic function.  

\begin{thm}\label{lower and upper bounds}
The asymptotic behavior of $N_n(r)$ is described by the function  
\[
c(r) = \lim_{n\to\infty} \frac{N_n(r)}{n},
\]
whose lower and upper bounds (and known exact values) are depicted in Figure~\ref{fig:graphupdate}.
\end{thm}

The proof of Theorem~\ref{lower and upper bounds} draws upon a collection of examples and short arguments, all presented in Section~\ref{sec: others}. There, we include a summary Table~\ref{tab:results} (including Theorem \ref{thm:half} and Theorem \ref{thm:fourth}) and the corresponding graph illustrating our current partial knowledge of $c(r)$. Together, these results provide a comprehensive picture of the currently known asymptotic behavior of $c(r)$ and highlight the open directions for further study.

We now outline the rest of the paper.  
\\In \cref{sec:half} we prove Theorem \ref{thm:half}, \\in \cref{sec:fourth} we prove Theorem \ref{thm:fourth}, and \\in \cref{sec: others} we prove Theorem \ref{lower and upper bounds} and include questions and remarks for future work.

\section{Proof of Theorem \ref{thm:half}}\label{sec:half}

\textit{Proof of Theorem \ref{thm:half}.} To prove that $N_{n}(\frac{1}{2})=\lceil \frac{n}{3}\rceil +1$, we consider the three remainder classes of  $n$ modulo $3$ separately. We set up a table to understand what needs to be proved for each remainder class without using the floor and the ceiling notation.

\begin{table}[ht]
\centering 
\begin{tabular}{|c|c|c|} 
\hline

Three cases according to & Prove that this many  & Find a point set where\\
 mod 3 of the cardinality& points can be covered &   no circle of radius $\frac 12$ can \\ 
  of the given point set& by a circle of radius $\frac 12$ &  cover this many points \\ [0.5ex]
\hline \hline 
$3m$ 		&Claim 1  \quad $m+1$ 	& Example 1 \quad $m+2$  \\ [0.5ex]
\hline
$3m+ 1$ 	& Claim 2 \quad $m+2$ 	& Example 1b \quad $m+3$\\ [0.5ex]
\hline
$3m +2$	 & Claim 3 \quad $m+2$ 	&   Example 1c \quad $m+3$  \\  [0.5ex]
\hline 
\end{tabular}
\caption{Restating Theorem \ref{thm:half}}
\label{table:nonlin} 
\end{table}

It suffices to verify  Claims 1-3 and to find  Examples 1, 1b, 1c of the above table. In fact, we need to prove only Claim 2 and present Example 1. Indeed, on the one hand Claim 2 implies Claim 1 by adding to the $3m$-point set one random extra point, on the other hand, it also implies Claim 3 as the latter one is a weaker claim. Examples 1b and 1c  will be simple modifications of Example 1.

\textbf{Proof of Claim 2.} Let $P$ be a set of $3m+1$ points whose diameter is at most $1$. Let $S$ with center $O$ be the smallest circle containing $P$ and let $H$ be the convex hull of those points which are on the boundary of $S$. According to Jung’s theorem, the radius $r$ of the circle $S$ is at most $\tfrac 1 {\sqrt3} \approx 0.577$. We may assume $r > \tfrac 12$, otherwise  all $3m+1$ points can be covered with a circle of radius $\tfrac 12$.

\medskip

\begin{wrapfigure}{r}{0.4 \textwidth}
\centering
\includegraphics[scale=.45]{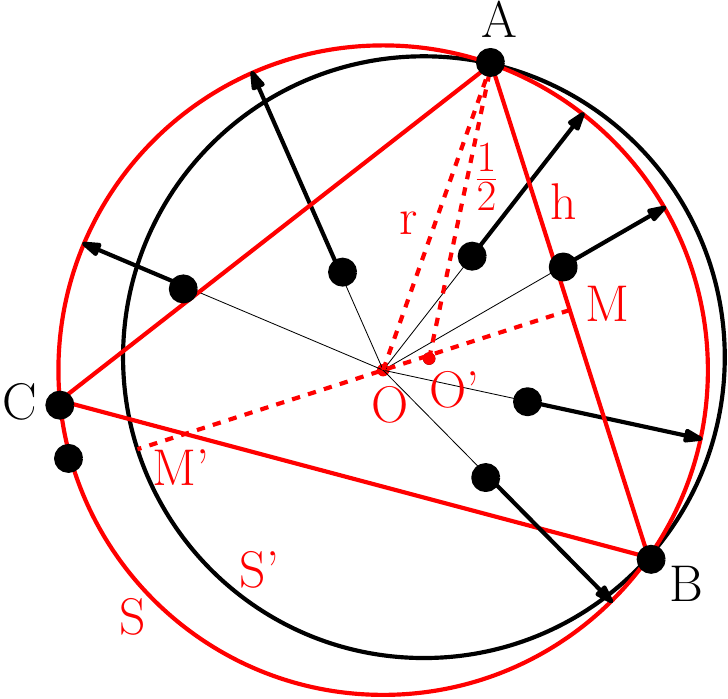}
\captionsetup{font=small, width=0.35\textwidth}
\caption{Proof of Claim 2.}
\label{case3n+1}
\end{wrapfigure}

Since the distance between any two points of $P$ is at most $1$, no two points of $P$ can be diagonally opposite points of $S$.  Thus, the center $O$ cannot be on the boundary of the polygon $H$.  If the center $O$ is an exterior point of $H$, then $S$ is not the smallest circle containing $P$, a contradiction. 
 Finally we conclude that, in any triangulation of the polygon $H$, there is a triangle whose interior contains $O$ and thus this triangle is  an acute triangle. Let us label the vertices of this acute triangle in clockwise order by   $A$, $B$ and $C$.  The significance of triangle $ABC$ being acute is that the shorter closed circular arcs $AB$, $BC$ and $CA$ of $S$ cover the circle $S$. 

 Finally, we explain which circle of radius $\tfrac 12$ will contain $m+2$ points. Let us project radially from $O$ all points of $P$ to $S$. If a point coincides with $O$, then we are free to project it to anywhere on $S$. The points $A$, $B$ and $C$ stay fixed under this projection. By the pigeonhole principle, one of the shorter closed circular arcs $AB$, $BC$ and $CA$ contains $n+2$ points. Assume arc $AB$ is this arc. For the proof of Claim 2 it is enough to verify that all the preimages of the projected points of arc $AB$ are covered by the  circle $S'$ of radius $\frac12$ with center $O'$ in $S$ and  passing through $A$ and $B$.
 
 Covering the preimages will readily follow from (1) $S'$ contains the shorter arc $AB$ of  circle $S$ and (2) $S'$ contains the center $O$ of $S$  (see Fig \ref{case3n+1}).
 (1) holds since $r > \frac 12$.
 To prove (2) let  $M$ be the midpoint of the chord $AB$ and let $M'$ be the diagonally opposite point of $M$ on circle $S'$. Center $O'$ belongs to the segment $OM$. Thus, all we need to verify is $\frac 12  \geq |OO'|$. Let us introduce $h$ as the length of segment $AM$.  For a fixed $h$, $OO'$ is largest when $ r =\frac 1 {\sqrt3}$. Thus, it is enough to verify  that 
 $\frac{1}{2} \geq \sqrt{\frac{1}{3}-h^2} - \sqrt {\frac{1}{4}-h^2}$
  for all $h \in (0,  \frac{1}{2} ]$.  This inequality simplifies to $\sqrt {\frac{1}{4} - h^2}  \geq - \frac{1}{6}$, which holds for $h \in  (0, \frac{1}{2}]$. \claimqed

\bigskip

Now, we turn to proving that $\lceil \frac{n}{3}\rceil +1$ is not only a lower bound, but it is also an upper bound for $N_{n}(\frac{1}{2})$. This requires constructing  Example 1 in Table \ref {table:nonlin}.

\medskip

\begin{tcolorbox}[colback=white, colframe=black, boxrule=.8pt,
  breakable, enhanced jigsaw]

% --- Top: full-width intro paragraph ---
\textit{We construct a diameter 1 point-set of $3m$ points, where a circle of
radius $\frac12$ cannot cover $m+2$ points; in other words,}

\begin{ex}
\textbf{We have that $N_{3m}(\frac{1}{2}) < m+2$.}
\label{ex:Emplus2}
\end{ex}

% \bigskip

% --- Next: paragraph + first figure side-by-side ---
\begin{tcolorbox}[colback=white, colframe=white, boxrule=0pt,
  enhanced jigsaw, sidebyside, righthand width=0.35\textwidth,
  sidebyside align=top seam, middle=0pt]

Consider the following $3m$-point arrangement. Take the Reuleaux triangle
$UVW$ with diameter $1$. On the boundary of this triangle, outside of the
centrally placed circle of radius $\frac12$, we will place a total of $3m$
points, $m$ on each side, so that the final arrangement has  3-fold
symmetry. Note that the final arrangement will not have the $m$ points on
each side equally spaced (Figure \ref{fig:mplus2}). Since the final
arrangement will have 3-fold symmetry, by describing the points along side
$UV$, this determines the placement of the points on sides $VW$ and $UW$.
The following is the recursive instruction to place $m$ points
$X_1,\dots,X_m$ (in increasing order of their distance to $V$) on the curved side $UV$. Let $Y_1,\dots,Y_m$ and
$Z_1,\dots,Z_m$ be the points on $VW$ and $UW$, respectively.

\tcblower
\centering
\includegraphics[scale=0.5]{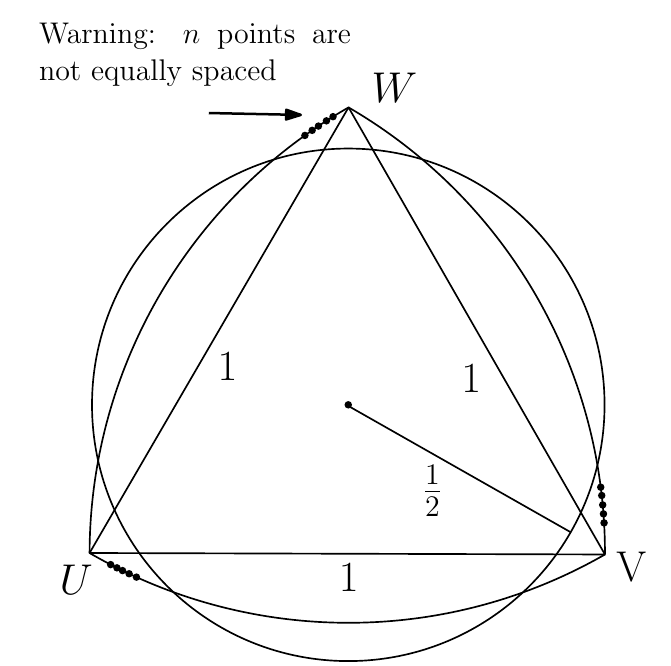}
\captionsetup{font=small, width=0.9\textwidth}
\captionof{figure}{Example \ref{ex:Emplus2}}
\label{fig:mplus2}

\end{tcolorbox}

% --- Middle: itemized list full width ---
\bigskip
\begin{itemize}[leftmargin=1em]
  \item First choose $X_1$ to be very close to $U$. Note that the choice of $X_1$ will also determine $Y_1$ and $Z_1$. As part of the iterative process, when we choose sequential $X_i$’s, we revisit the choice of the initial closeness and make necessary adjustments (this adjustment will be explained in the next paragraph).
  \item Now the recursive step: once we have $X_i$, choose $X_{i+1}$ (which is the same as choosing $Y_{i+1}$) so that the circumradius of the triangle $X_iY_iY_{i+1}$ is greater than $\tfrac12$ (see Figure \ref{fig:detailsmplus2}). Can we really choose $X_{i+1}$ (and thus $Y_{i+1}$) like this? Answer is yes, assuming that $x_i$ is sufficiently close to $U$ and thus $Y_i$ is automatically at the same sufficiently close distance from $V$. The exact reason is the following: First notice that no matter where we choose $Y_{i+1}$, the  triangle $UVY_{i+1}$ is isosceles with $UV = UY_{i+1} =1$, therefore  triangle $UVY_{i+1}$ always has a circumradius $> \frac 12$. Now, if we relocate  $U$ and $V$  to $X_i$ and $X_{i+1}$ within a sufficiently close neighbourhood of $U$ and $V$ the circumradius of the triangle $X_iY_iY_{i+1}$ remains $>\frac 12$. Finally it is very important to guaranty for the next iteration that  $Y_{i+1}$ is as close to $V$ as we wish even after the $i$-th iterations. To ensure this, before doing the next iteration, we might need to go back to the beginning and  choose $X_1$ much much  closer to $U$ and repeat the entire process. This is the adjustment which was mentioned at explaining the initial step of the iterative process.
\end{itemize}

% --- Bottom: second figure full width ---
\bigskip
\begin{center}
\includegraphics[scale=0.6]{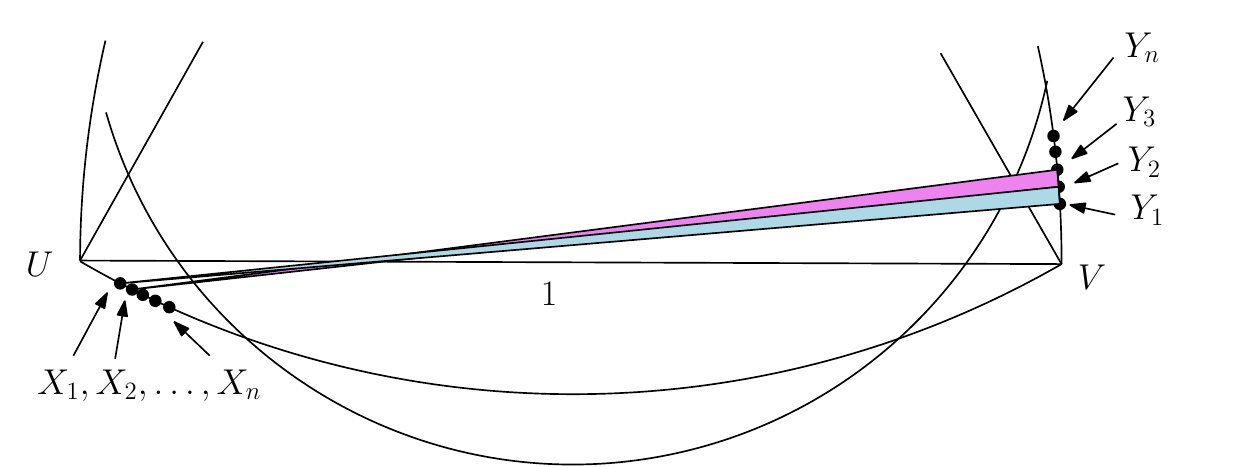}
\captionsetup{font=small}
\captionof{figure}{Enlarged lower part of Figure \ref{fig:mplus2}}
\label{fig:detailsmplus2}
\end{center}

\bigskip
Now we must show that  no circle of radius $\tfrac{1}{2}$ can cover more than $m+1$ points of the just constructed $3m$  points. Every triangle $X_iY_jZ_k$ has a circumradius  $>\frac12$, thus it is enough to show that no circle of radius $\tfrac{1}{2}$ can cover more than $m+1$  points from two of the sides of the triangle $UVW$,  say from $UV$ and $VW$. Otherwise, there are two indices $i<j$ with the property that four points $X_i,Y_i,X_j,Y_j$ are covered by the circle of radius $\frac12$. This implies that  the triangle $X_i,Y_i,Y_{i+1}$ is also covered by this circle, which is impossible since by construction this triangle has a circumradius greater than $\tfrac12$.  \claimqed

\end{tcolorbox}

\medskip

\section{Proof of Theorem \ref{thm:fourth}}\label{sec:fourth}

Before giving the proof, we make a comment on work of Bateman and Erd\H{o}s \cite{BE51}, who proved that $N_7(\frac{1}{4})= 2$.  Certainly, $N_7(\frac{1}{4})\leq 2$ is true by placing one point at each of the 7 vertices of a regular heptagon with diameter $1$. Hence the challenge lies in showing $N_7(\frac{1}{4})\geq 2$, which is equivalent to: Show that there does not exist a $7$-point set $P$ with the property that $\frac{1}{2}<|xy|\leq 1$ for every distinct $x,y\in P$. This is the exact version which was solved by Bateman and Erd\H{o}s.

\textit{Proof of Theorem \ref{thm:fourth}.} Due to the above result of Bateman and Erd\H{o}s  we can assume $n\neq 7$. As it was mentioned already in Section \ref{sec:intro}, $\lceil \frac{n}{7} \rceil$ is a lower bound for $N_n(\frac{1}{4})$, in fact this follows immediately from the circle covering depicted in the lower left corner of Figure \ref{Yanchao}. Indeed,  $7$ congruent circles of radii $\frac 14$ cover a regular hexagon whose opposite sides are at distance 1. This implies, by the pigeonhole principle that $N_n(r) \geq \lceil \frac n7 \rceil$, when $r  \geq \frac{1}{4}$.

Concerning upper bounds for  $N_n(\frac{1}{4})$, we provide examples in which $\lceil \frac{n}{7} \rceil+1$ points cannot be covered by a circle of radius $\frac{1}{4}$ when $n$ is not a multiple of $7$ and $\lceil \frac{n}{7} \rceil+2$ points cannot be covered when $n$ is a multiple of $7$. We distinguish three cases, depending on $n$, for which we find different examples:

\textbf{Case 1: $n\leq 5$}.  
\begin{tcolorbox}[colback=white, colframe=black, boxrule=0.8pt, breakable, enhanced jigsaw] 
{\it For each $n \leq 5$, we need $n$ points of diameter at most $1$ so that $ \lceil \frac{n}{7} \rceil+1= 2$ points cannot be covered with a circle of radius $\frac14$.} We can see that a regular $n$-gon of diameter $1$ has side lengths greater than $\frac{1}{2}$, and hence placing a point on each vertex gives $N_n(\frac{1}{4}) = 1$. 
\end{tcolorbox}

\textbf{Case 2: $n = 6$}. 
\begin{tcolorbox}[colback=white, colframe=black, boxrule=0.8pt, breakable, enhanced jigsaw]

{\it We need 6 points of diameter at most $1$ so that $ \lceil \frac{6}{7} \rceil+1= 2$ points cannot be covered with a circle of radius $\frac14$.} We can take a circle of radius $\frac{1}{2}+\epsilon$ for some small $\epsilon>0$, and place $5$ points uniformly around the circle with one point in the center. Since no two points on the circle are antipodal, with a sufficiently small $\epsilon$, these $6$ points have diameter at most $1$. Similar to the regular pentagon, any circle of radius $\frac{1}{4}$ can only cover one point. 
\end{tcolorbox}

\textbf{Case 3: $n>7$.}

\begin{tcolorbox}[colback=white, colframe=black, boxrule=0.8pt, breakable, enhanced jigsaw]

\begin{ex}
For each $n>7$, we describe a point set of $n$ points with diameter $<1$ so that analysis of the point arrangement reveals that a circle of
radius $\frac14$ cannot cover

\medskip

A) $\lceil \frac{n}{7}\rceil+1$ points when $n$ is not a multiple of $7$ and 

B) $\lceil\frac{n}{7}\rceil+2$ points when $n$ is a multiple of $7$. \label{ex:fourth}
\end{ex}

\medskip

Let $S$ be the following $n$-point set. First, distribute $\lceil \tfrac{n}{7}\rceil$ points uniformly along the boundary of a small circle with radius $\epsilon>0$ (see Figure~\ref{rfourth}). Then distribute uniformly  the remaining points along the boundary of a circle of radius $\tfrac{1}{2}$.
\medskip

Quantitatively, this means that if $n= 7k+m$ (for $m\in\{1, ...,6\}$) then $\lceil \tfrac{n}{7}\rceil = k+1$, the $\epsilon$-circle contains $k+1$ points, the outer circle contains $6k+m-1$ points. Also notice that a circle of radius $\frac 14$ can cover at most $\lfloor \frac{6k+m-1}6\rfloor +1 = k+1$ points on the outer circle. 

\medskip
A similar chain of thought means that if $n=7k$, then $\lceil \frac{n}{7}\rceil = k$, the $\epsilon$-circle contains $k$ points, the outer circle contains $6k$ points. Also notice that a circle of radius $\frac 14$ can cover at most $ \frac{6k}6 +1 = k+1$ points on the outer circle.

\medskip
\begin{center}
\includegraphics[scale=1.4]{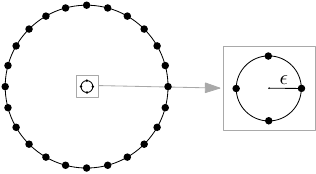}
\captionsetup{font=small}
\captionof{figure}{Example \ref{ex:fourth}}
\label{rfourth}
\end{center}

\medskip

We will finish the analysis of this case by looking at the following three subcases:

\medskip    

\begin{enumerate}[leftmargin=1.5em]
  \item[(i)] \textit{Circle $O$ covers at least two points of $S$ from the outer circle.} In this case, when $\epsilon$ is sufficiently small, circle $O$ cannot cover any point of the inner circle. This implies that $O$ can cover at most  $k+1$ points, regardless if $n$ is divisible by 7 or not. Note that $k+1$ is  $\lceil \tfrac{n}{7}\rceil$ if $n=7k+m$, i.e. when $n$ is not divisible by $7$ and it is  $\lceil \tfrac{n}{7}\rceil +1$ if $n=7k$, i.e. when $n$ is divisible by $7$.

  \item[(ii)] \textit{Circle $O$ doesn't cover any point from the outer circle.} In this case circle $O$ covers at most $\left\lceil \tfrac{n}{7} \right\rceil$ points of $S$, again regardless if $n$ is divisible by $7$ or not.

  \item[(iii)] \textit{Circle $O$ covers exactly $1$ point of $S$ from the outer circle.} Since $n>7$, there are at least two points of $S$ on the inner circle of radius $\epsilon$. Since the points of $S$ are uniformly distributed along the inner circle, we have at least one point from the inner circle is not covered by $O$. Thus, if $n$ is not divisible by $7$ and has the form $n=7k +m$, then at most $k+1$ (which is $\lceil \tfrac{n}{7}\rceil$) points can be covered by a circle of radius $\frac14$. On the other hand, if $n$ is divisible by 7 and has the form $n=7k$ then at most $k$ points (which is $\lceil \tfrac{n}{7}\rceil$) can be covered by a circle of radius $\frac14$.
  
\end{enumerate}
\claimqed

\end{tcolorbox}

\section{Proof of Theorem \ref{lower and upper bounds}}\label{sec: others}

\begin{table}[h!]
    \centering
    \renewcommand{\arraystretch}{1.65} % increases the height of each row
    \setlength{\tabcolsep}{2.75pt} % sets the column separation
    
    \begin{threeparttable}
        \begin{tabular}{|l|l|l|l|}
            \hline
            Graph pieces & Over the intervals of $r$ & $N_n(r)$  & $c(r)$\\
            \hline \hline
             1.  Segment \textbf{A} (Jung)& $r \in [\frac{1}{\sqrt{3}},1] \approx [.577, 1]$ & $N_n(r) = n$ & $c(r)=1$\\
            \hline
            2.    Segment \textbf{B}  (Prop \ref{thm: B}) & $r \in [\frac{\sqrt{3}}{4}, \frac{1}{2})  \approx [.433, .5]$ & $N_n(r) = \lceil \frac{n}{3} \rceil$ & $c(r)=\frac13$\\
            \hline
            3.   Point \textbf{b}  (Thm \ref{thm:half})  & $r= \frac{1}{2}$  & $N_{n}(\frac{1}{2})= \lceil\frac{n}{3}\rceil+1$ & $c(\frac12 )=\frac13$\\
            \hline 
            4.   Point \textbf{c} (Prop \ref{thm: c})  & $r= \frac{\sqrt{2}}{4} \approx .354$ & $N_{4n}(\frac{\sqrt{2}}{4})\in  \{n, n+1 \}$ & $c(\frac{\sqrt{2}}{4})= \frac 14$\\
            \hline
            5.   Point \textbf{d} (Thm \ref{thm:fourth})  & $r= \frac{1}{4}$ & $N_n(\frac{1}{4}) = \lceil \frac{n}{7} \rceil, N_7(\frac{1}{4})=2$ &$c(\frac{1}{4})=\frac17$\\
            \hline
            6.   Segment \textbf{E} (Prop. \ref{thm: E})  & $r \in [\frac{\sqrt{3}-1}{\sqrt{2}},\frac{1}{\sqrt{3}}) \approx [.518,.577]$ & $N_n(r) \leq \lceil \frac{2n}{3} \rceil$ &$c(r) \leq \frac23$\\
            \hline
            7.    Curve \textbf{F} (Prop. \ref{thm: F}) & $r \in (0,\frac{1}{2})$  & $N_n(r) \leq \lfloor \frac{n}{\pi}\sin^{-1}(2r)+1 \rfloor$ & see \tnote{*}\\
            \hline
            8.    Segment \textbf{G} (Prop. \ref{thm: G})& $r \in [\frac{1}{2},\frac{\sqrt{3}-1}{\sqrt{2}}) \approx [.5,.518]$ & $N_n(r) \leq \lceil \frac{3n}{5} \rceil$ & $c(r) \leq \frac 35$\\
            \hline
            9.   Segment \textbf{H} (Prop. \ref{thm: H}) & $r \in [\sqrt{\frac{13}{48}},\frac{1}{\sqrt{3}}]\approx [.52,.577]$ & $N_n(r) \geq \lceil \frac{n}{2} \rceil$& $c(r) \geq \frac 12$\\
            \hline
            10.   Function \textbf{I} (Fig. \ref{fig:graphupdate}) & $r \in  [.209,.433]$ & see \tnote{**} & see \tnote{**}\\
            \hline
        \end{tabular}
        
        \begin{tablenotes}
            \item[*] $c(r) \leq  \frac{\sin^{-1}(2r)} {\pi}$
            \hskip 1in    ** The gray function \textbf{I} in  Fig. \ref{fig:graphupdate} was derived in Sec. \ref{sec:intro}.
        \end{tablenotes}
        \bigskip
        
    \caption{A summary of the results concerning  $N_n(r)$ and of $c(r)$ }
    \label{tab:results}
    \end{threeparttable}
    
\end{table}

\begin{figure}[h!]
\centering
\includegraphics[scale = 1.1]{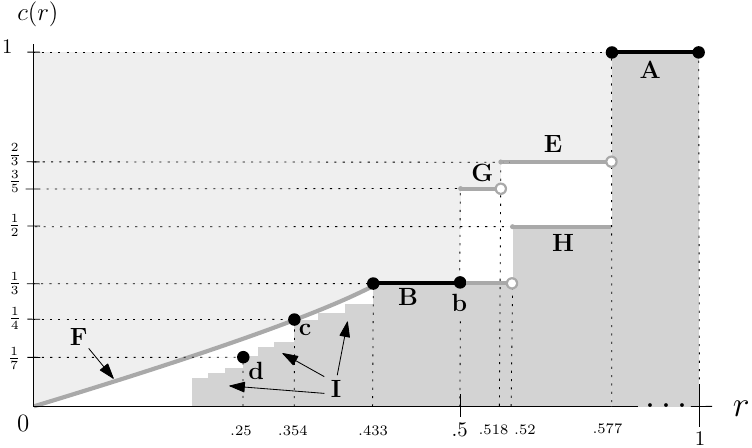}
    \caption{Partial knowledge of the graph of  $c(r) := \lim_{n\to \infty} \frac{N_n(r)}{n}$}
\label{fig:graphupdate}
\end{figure}

Table \ref{tab:results} contains what we know about the graphs of the functions $N_n(r)$ and $c(r)$, where $c(r)$ is defined as $c(r) := \lim_{n\to \infty} \frac{N_n(r)}{n}$. Figure \ref{fig:graphupdate} helps to visualize the same information graphically.  To help the reader navigate through the theorems and examples in this paper, certain intervals and points in Figure \ref{fig:graphupdate} will be labeled with bold letters from {\bf A} to {\bf I} and these letters will be used as reference letters in the theorems throughout the paper (uppercase for intervals and lowercase for points). Exact values of the function $c(r)$ are indicated in black, lower bounds in dark gray, and upper bounds in light gray. (Note: the $r$-axis is not drawn to scale for visual interest). 

Next we briefly interpret row by row the content of Table \ref{tab:results}. 
Rows 1, 3 and 5 are the same as Jung's theorem (Theorem \ref{jung}) and Theorems \ref{thm:half} and \ref{thm:fourth}, which were already proved in this paper. Row 10 is the step function derived in \cref{sec:intro} from the work of Liu \cite{Y22} in \cref{Yanchao}. This explain why segments A,  point b, and point d are part of the graph of the function $c(r)$. It turns out that segment B and point c are also points of the graph of $c(r)$. In fact, we have the following.

\begin{prop}\label{thm: B} {\bf [Interval \hyperref[fig:graphupdate]{\textbf{B}}]} Row $2$ of Table \ref{tab:results} says that the open segment B is part of the graph of $c(r)$ due to $N_n(r) = \lceil \frac n3 \rceil$ for any $r \in [ \frac{\sqrt3}4, \frac 12 )$.
\end{prop}

\begin{tcolorbox}[colback=white, colframe=white, boxrule=0pt,
  enhanced jigsaw, sidebyside, righthand width=0.3\textwidth,
  sidebyside align=top seam, lefthand ratio=0.65, middle=0pt]

{\it Proof.} Showing $\lceil \frac n3 \rceil\leq N_n(r)$ follows from the idea of  Borsuk \cite{B33,Z21} who proved that any set $C$ of diameter $1$ can be partitioned into three pieces each with diameter $< 1$ . Borsuk's  enclosed $C$ in a regular hexagon whose opposite sides are at a distance $1$, then  partitioned this hexagon into three congruent pentagons by dropping perpendiculars from every second side (Figure \ref{intervalB}). By pigeonhole principle, one of the pentagons contains at least $\lceil \frac{n}{3}\rceil$ points. Since the circumradii of these pentagons are equal to $  \frac {\sqrt 3}4$ we obtain $N_n(r) \geq \lceil \frac n3 \rceil$, when $r \geq \frac {\sqrt 3}4$. Finally, the balanced distribution of $n$ points among the three vertices of an equilateral triangle shows that  $\lceil \frac n3 \rceil$ is also an upper bound for $N_n(r)$. \hfill \claimqed

\tcblower
\centering
\includegraphics[scale=0.7]{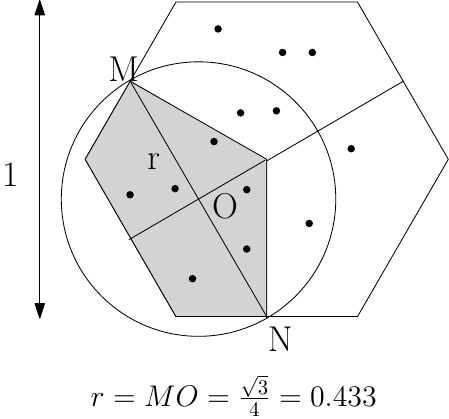}
\captionsetup{font=small, width=0.9\textwidth}
\captionof{figure}{Proposition \ref{thm: B}}
\label{intervalB}

\end{tcolorbox}

\medskip

\begin{prop}\label{thm: c}{\bf[Point \hyperref[fig:graphupdate]{c}]} Row $4$ of Table \ref{tab:results} says that point c is on the graph of c(r) due to 
$N_n(\frac{\sqrt{2}}{4}) = \lceil \frac{n}{4}\rceil$ for $n\not \equiv 0\pmod{4}$ and $N_{4n}(\frac{\sqrt{2}}{4}) \in \{n,n+1\}$. 
\end{prop}

\textit{Proof.} Let us enclose the point set of diameter $1$ in a square of side length $1$. Partition the square into four squares of side length $\frac 12$  
 
By pigeonhole principle, one of them contains at least $\frac 14  n $ points.  Thus $N_n(r) \geq \lceil \frac n4 \rceil$ whenever $ r \geq \frac {\sqrt2}4$. Finally,  distribute uniformly $n$ points on a circle of radius $r = \frac{\sqrt{2}}{4}$. It is easy to see that a circle of radius 
$\frac{\sqrt{2}}{4}$ can cover at most $\frac n4$ points if $n$ is a multiple of 4, otherwise it can cover at most $\lceil \frac{n}{4} \rceil$ as desired. 
\claimqed

\medskip

\begin{prop}\label{thm: E}{\bf[Segment \hyperref[fig:graphupdate]{E}]} Row 6 of Table \ref{tab:results} says that left-closed right-open segment E is an upper bound for the graph of $c(r)$.
\end{prop}

\begin{tcolorbox}[colback=white, colframe=black, boxrule=0.8pt, breakable, enhanced jigsaw]

\noindent\textit{We give a unit-diameter point set, where essentially}

\begin{itemize}[leftmargin=2em]
  \item[(a)] any circle of radius less than one-half cannot cover more than one-third of the points, 
  \item[(b)] any circle whose radius is at least one-half, but less than the Jung radius, cannot cover more than two-thirds of the points.
\end{itemize}
  Take an equilateral triangle of side length $1$ with 
the balanced distribution of  $n$ points among the three vertices of this triangle.  Certainly, this point set has diameter $1$, and if $r \in (0,\frac{1}{2})$ then a circle of radius $r$ can cover at most one vertex. Likewise, if $r \in [\frac{1}{2}, \frac{1}{\sqrt{3}})$ then a circle of radius $r$ can cover at most two vertices. In the first case no circle can cover more than $\lceil \frac{n}{3} \rceil$ points, and in the second, no circle can cover more than $\lceil \frac{2n}{3} \rceil$ points. \claimqed
\end{tcolorbox}

\begin{prop}\label{thm: F}{\bf[Curve \hyperref[fig:graphupdate]{F}]} Row 7 of Table \ref{tab:results} gives a better upper bound for $c(r)$ over the interval $(0, \frac 12)$, than the one that was used in Section \ref{sec:intro}.
\end{prop}

\begin{tcolorbox}[colback=white, colframe=black, boxrule=0.8pt, breakable, enhanced jigsaw]

\textit{In view of Theorem~\ref{thm:half}, we already know that any circle of radius less than one-half cannot cover more than one-third of the points. It is expected that for small radii, the number of points such a circle can cover is much less than one-third of the points. We give a unit-diameter set where indeed far fewer points can be covered. More precisely, we show that $N_n(r) \leq \lfloor \tfrac{n}{\pi}\sin^{-1}(2r) \rfloor + 1$ for any $0 < r < \tfrac{1}{2}$.}

\begin{tcolorbox}[colback=white, colframe=white, boxrule=0pt,
  enhanced jigsaw, sidebyside, righthand width=0.35\textwidth,
  sidebyside align=top seam, lefthand ratio=0.65, middle=0pt]

\textit{Proof.}
Consider a circle $C$ with radius $\tfrac{1}{2}$ and uniformly distribute all $n$ points around the boundary of $C$ (see Figure~\ref{fig:exampleF}). Certainly this point set has diameter at most $1$. 

\medskip

If $r \in (0, \tfrac{1}{2})$, then a circle of radius $r$, say $C'$, can cover an arc of length at most $\sin^{-1}(2r)$. This maximum is achieved when the intersection of $C$ and $C'$ occurs at antipodal points of $C'$. Since the point set is equally spaced around $C$, $C'$ covers at most $\lfloor \tfrac{n}{\pi}\sin^{-1}(2r) \rfloor + 1$ points. \claimqed

\tcblower
\centering
\includegraphics[scale=0.6]{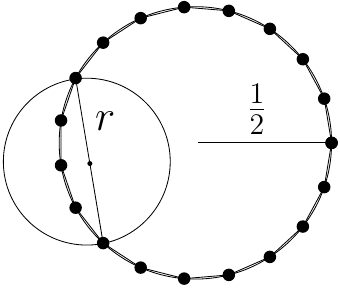}
\captionsetup{font=small, width=0.95\textwidth}
\captionof{figure}{Example in Proposition~\ref{thm: F} with $n=19$}
\label{fig:exampleF}

\end{tcolorbox}

\end{tcolorbox}

\begin{prop}\label{thm: G}{\bf[Segment \hyperref[fig:graphupdate]{G}]} Row 8 of Table \ref{tab:results} gives a better upper bound for $c(r)$ over the interval $(0, \frac 12)$, than the one which was used in Section \ref{sec:intro}.
\end{prop}

\begin{tcolorbox}[colback=white, colframe=black, boxrule=0.8pt, breakable, enhanced jigsaw]

\noindent\textit{We give a unit-diameter point set such that essentially any circle of radius less than $\frac{\sqrt{3}-1}{\sqrt{2}}\approx 0.518$ cannot cover more than three-fifths of the points. More precisely, we show that $N_n(r) \leq \lceil \frac{3n}{5} \rceil$ for any $r \in [\frac{1}{2},\frac{\sqrt{3}-1}{\sqrt{2}})$.   }

Let $r < \frac{\sqrt{3}-1}{\sqrt{2}}$ and let $O$ be a circle of radius $r$. Consider the Reuleaux triangle $UVW$ of diameter $1$. Let $M_{1}$ and $M_{2}$ be the midpoints of the circular arcs $UV$ and $UW$ respectively. Now consider the $5$ point arrangement $P$ by placing a point on the vertices $\{U,V,W,M_1,M_2\}$. We now take an $n$-point balanced arrangement of $P$ where if $n$ is not a multiple of $5$, we place the remaining points each on a distinct vertex, favoring $\{U,V,W\}$. Since $r < \frac{\sqrt{3}-1}{\sqrt{2}} < \frac{1}{\sqrt{3}}$ as mentioned in the Example from Proposition \ref{thm: E}, $O$ covers at most two of the vertices in $\{U,V,W\}$. Now suppose for contradiction $O$ can cover more than $\lceil \frac{3n}{5}\rceil$ of the points. Then without loss of generality $O$ must cover $\{U,V,M_1,M_2\}$. It is an easy computation to verify that this would require $r\geq \frac{\sqrt{3}-1}{\sqrt{2}}$, a contradiction. \claimqed
\end{tcolorbox}

\begin{prop}\label{thm: H}{\bf[Segment \hyperref[fig:graphupdate]{H}]} Row 9 of Table \ref{tab:results} says that the closed horizontal segment H is a lower bound of $c(r)$ over its interval due to $N_n(r) \geq \lceil \frac n2 \rceil$ for any $r\in [\sqrt{\frac{13}{48}}, \frac{\sqrt{3}}{3})$..
\end{prop}

\begin{tcolorbox}[colback=white, colframe=white, boxrule=0pt,
  enhanced jigsaw, sidebyside, righthand width=0.35\textwidth,
  sidebyside align=top seam, lefthand ratio=0.65, middle=0pt]

{\it Proof.} All what we need is to remember that according to Borsuk \cite{B33, Z21}, the regular hexagon with opposite sides at distance $1$ covers every unit diameter set. Let us cut the hexagon into two congruent pentagons by the perpendicular bisector of one of its sides. Figure \ref{felBorsuk} shows that the circumradius of this pentagon is  $r = \sqrt \frac {13}{48} \approx 0.52 $. By the pigeonhole principle, one of the pentagons contains at least half of the points, thus $N_n(r) \geq \lceil \frac n2 \rceil$. \claimqed

\tcblower
\centering
\includegraphics[scale=0.3]{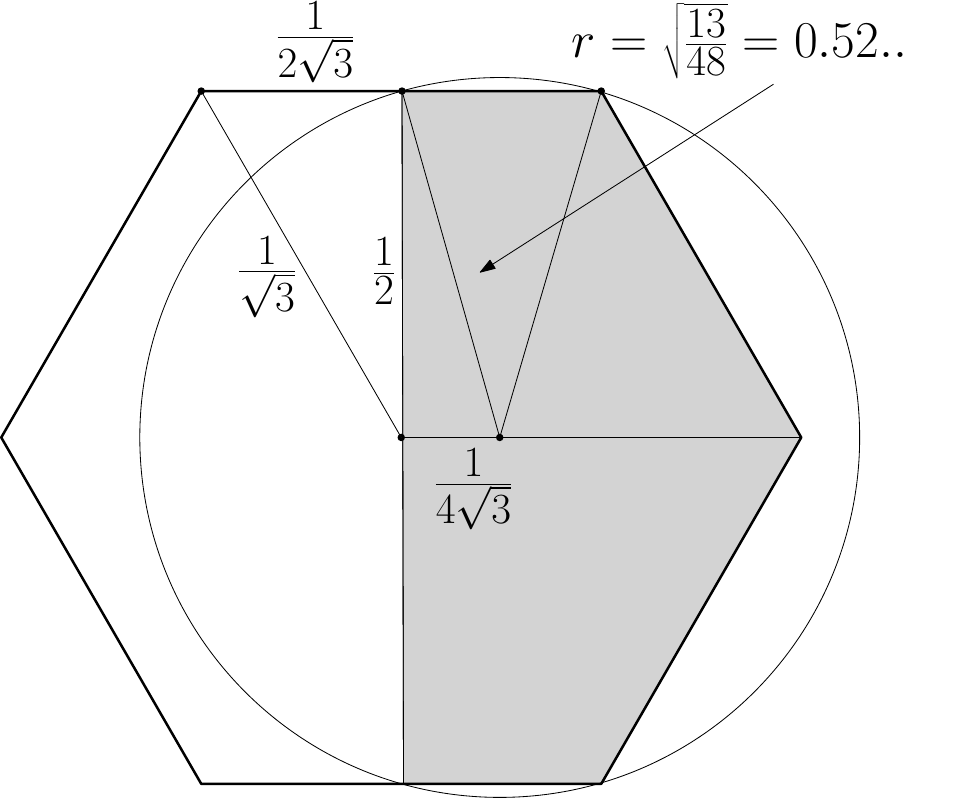}
\captionsetup{font=small, width=0.9\textwidth}
\captionof{figure}{Proposition~\ref{thm: H} }
\label{felBorsuk}

\end{tcolorbox}

\section{Open questions and future work} \label{future directions}

1. Finding additional radii $r$, for which the exact function values of $N_n(r)$ can be determined remains the most interesting question. 

2. In view of Theorem \ref{thm:fourth}, $r= \frac 14$ could be the next such radius, as at this point we know all values of $N_n(\frac14)$ except for those $n$ which are multiples of $7$. We proved that $N_{7k}(\frac 14) = k$ or $k+1$. Bateman and Erdos handled the first  case by proving that  $N_7(\frac 14)=2$.  Generalizing their result  for other multiples of $n$ seems to be an important, but a formidable task.

3. Jung's theorem, which inspired this paper, was generalized to  higher-dimensional spaces.  The general theorem for dimension d states that every unit diameter point set can be covered by a closed ball of  radius \(\leq \sqrt{\frac{d}{2(d+1)}}\). It is also natural to study maximal partial coverages with smaller spheres.

4. The lower bounds for $N_n(r)$ were using circle coverings of certain universal covers of the unit diameter sets.
One might want to identify and study collections of convex shapes (called {\it covering collections}) such that any set of unit diameter  can be covered by at least one shape in the collection. For our purposes those covering collection are useful when we can prove statements like ``each shape in the collection can be covered by $k$ congruent circles of a certain radius $r$". If, for some $k$, the radius $r$ is smaller than that in case of a single universal cover, then this yields a stronger bound on $N_n(r)$.

5. Our partial covering problem has an analogous variant for every pair of convex sets $(K, L)$. Instead of unit diameter point sets, consider sets of $n$ points in $K$. For a specific  $\lambda$ we would like to find the largest $k$ such that a homothetic copy $\lambda L$  can cover  $k$ points in every $n$-point set of $K$. It would be interesting to find  a pair $(K,L)$  which leads to elegant arguments and results.

\textbf{Data Availability.} Out manuscript has no associated data.

\end{document}